\documentclass{ifacconf}

\usepackage{graphicx}      % include this line if your document contains figures
\usepackage{natbib}        % required for bibliography

\usepackage{amsmath, amssymb, amsfonts}
\usepackage{adjustbox}
\usepackage{array}
\usepackage{etoolbox}     % for \patchcmd, used by the section-spacing knob below

\usepackage{xcolor}

\makeatletter
\def\lb@subsubsection{\thesubsubsection.\half@em}
\def\lb@empty@subsubsection{\thesubsubsection}

\def\subsubsection{\@startsection{subsubsection}{3}{\z@}{0.3\@bls \@plus .2\@bls \@minus .1\@bls}{0.3\@bls}{\normalsize\itshape}}

\newdimen\secheadskip
\patchcmd{\section}{6\p@ \@plus 0.5\@bls \@minus 4\p@}{\secheadskip \@plus 0.5\@bls \@minus 4\p@}{}{\ClassWarning{main_v02}{section before-skip patch FAILED}}
\patchcmd{\section}{6\p@ \@plus 1\p@ \@minus 2\p@}{\secheadskip \@plus 1\p@ \@minus 2\p@}{}{\ClassWarning{main_v02}{section after-skip patch FAILED}}
\patchcmd{\subsection}{6\p@ \@plus 0.5\@bls \@minus 4\p@}{\secheadskip \@plus 0.5\@bls \@minus 4\p@}{}{\ClassWarning{main_v02}{subsection before-skip patch FAILED}}
\patchcmd{\subsection}{6\p@ \@plus 1\p@ \@minus 2\p@}{\secheadskip \@plus 1\p@ \@minus 2\p@}{}{\ClassWarning{main_v02}{subsection after-skip patch FAILED}}
\newcounter{mythm}
\newcounter{myrem}
\renewenvironment{thm}[1][]{%
  \refstepcounter{mythm}%
  \if\relax\detokenize{#1}\relax
    \addcontentsline{toc}{section}{\protect\numberline{\themythm}Theorem}%
    \noindent\textbf{Theorem \themythm.}\ignorespaces
  \else
    \addcontentsline{toc}{section}{\protect\numberline{\themythm}Theorem (#1)}%
    \noindent\textit{Theorem \themythm} (#1).\ignorespaces
  \fi
}{\par}

\renewenvironment{rem}[1][]{%
  \refstepcounter{myrem}%
  \if\relax\detokenize{#1}\relax
    \noindent\textbf{Remark \themyrem.}\ignorespaces
  \else
    \noindent\textit{Remark \themyrem} (#1).\ignorespaces
  \fi
}{\par}
    
\makeatother

\begin{document}
\begin{frontmatter}

\title{Unifying Entropy Regularization in Optimal Control: \\
{\large From and Back to Classical Objectives via Iterated Soft Policies and Path Integral Solutions}\thanksref{footnoteinfo}} 

\thanks[footnoteinfo]{This work was supported by the Research Foundation Flanders (FWO) under
SBO grant no. S007723N}

\author[First,Second]{Ajinkya Bhole} 
\author[First,Second]{Mohammad Mahmoudi Filabadi} 
\author[First,Second]{Guillaume Crevecoeur}
\author[First,Second]{Tom Lefebvre}

\address[First]{Department of Electromechanical, Systems and Metal Engineering, Ghent University, Ghent, Belgium (e-mail: ajinkya.bhole@ugent.be).}
\address[Second]{Core lab MIRO, Flanders Make, Belgium.}

\begin{abstract}
This paper develops a unified perspective on several optimal control formulations through the lens of Kullback-Leibler (KL) regularization. We propose a central problem that separates the KL penalties on policies and transitions with independent weights, thus generalizing the standard trajectory-level KL-regularization used in probabilistic optimal control. This umbrella formulation recovers various control problems: the classical Stochastic Optimal Control (SOC), Risk-Sensitive Stochastic Optimal Control (RSOC), and their policy-based KL-regularized counterparts, termed soft-policy SOC and RSOC, which yield tractable surrogates. Beyond being regularized variants, these soft-policy formulations majorize the original SOC and RSOC, thus, iterating their solutions recovers the original objectives. We further identify a synchronized case of soft-policy RSOC where the policy and transition KL weights coincide, yielding a linear Bellman operator, path-integral solution, and compositionality---extending these computationally favourable properties to a broad class of control problems.
\end{abstract}


\end{frontmatter}
%===============================================================================

\section{Introduction}

Optimal control problems arise in a wide range of application areas, from robotics and autonomous systems (\cite{toussaint2009robot,williams2019modelTRO}) to finance and operations research (\cite{follmer2002convex}), where sequential optimal decisions must be made. They share a common structure: an agent shapes the (stochastic) behaviour of a dynamical system by minimizing a cost over a (finite) horizon. Such problems are typically addressed via dynamic programming, which recursively decomposes them into subproblems. These subproblems, however, remain highly nonlinear and exact closed-form solutions are limited to a few special cases (\cite{karny1996towards}). Motivated by this, decades of research have investigated alternative formulations whose solutions resemble those of classical optimal control while exhibiting more favourable computational properties.

\cite{karny1996towards} was among the first to address this, introducing the notion of probabilistic optimal control (\cite{karny2006fully}). Therein, optimal decision making is cast as a \emph{density-matching problem}, with the Kullback-Leibler (KL) divergence quantifying discrepancy between densities. Rather than minimizing a cost directly, the agent shapes the closed-loop trajectory density to match a desired target. The resulting solution mirrors dynamic programming, but optimization operators are replaced by expectations over known densities, \emph{thus improving the solution's tractability}.

K\'arn\`y's work is an early predecessor of the Control as Inference (CaI) paradigm~(\cite{toussaint2006probabilistic, toussaint2009robot, neumann2011variational, rawlik2012stochastic, levine2018reinforcement}), in which the cost is encoded into the system's probabilistic graphical model (PGM) via auxiliary \emph{optimality variables} whose observation probability is proportional to the exponential of the negative cost, making low-cost trajectories more probable. Two broad strategies have grown from this construction. The first treats policy derivation as a posterior inference problem on the PGM and, through exact message passing or variational inference over the joint trajectory distribution, reads the policy off the resulting posterior~(\cite{toussaint2006probabilistic, toussaint2009robot, levine2018reinforcement}). The second poses an explicit density-matching problem in which the policy is chosen to minimise a KL divergence between its induced trajectory distribution and the optimality posterior, with the direction of the divergence itself a deliberate modelling choice~(\cite{neumann2011variational, rawlik2012stochastic, lefebvre2024probabilistic, levine2013guided, levine2013variational}).

Both strategies share a conceptual reversal of classical optimal control design. In the classical viewpoint, a control objective is \emph{synthesised} first, and an optimal policy together with its closed-loop trajectory distribution emerges from minimising it. In CaI, the design instead begins from an inference or density-matching problem, and a policy is \emph{extracted} from its solution without ever writing down the classical control objective whose minimisation would yield that same policy. This reversal leaves a fundamental question unanswered: \emph{what interpretable, traditional control objective is actually being solved by these procedures?}

Partial answers exist within the density-matching family. \cite{neumann2011variational} first distinguished the two natural directions of this KL divergence, one in which the policy-induced distribution appears as the divergence's first argument (the so-called \emph{I-projection}) and one in which it appears as its second (the \emph{M-projection}), and observed that the two give rise to markedly different policies. Building on this, \cite{rawlik2012stochastic} showed that the I-projection yields an entropy-regularized variant of Stochastic Optimal Control (SOC) whose fixed-point iteration converges to the classical SOC policy, while \cite{lefebvre2024probabilistic} showed that the M-projection recovers exactly the policies of direct conditioning on optimality, with its iteration converging instead to the policy of Risk-Sensitive Stochastic Optimal Control (RSOC), wherein the cost enters exponentially to encode a risk-seeking attitude. Each step of this M-projection iteration also coincides with the M-step of an Expectation--Maximization update on the Maximum-Likelihood-of-optimality formulation of CaI~(\cite{toussaint2009robot,noorani2022probabilistic,watson2021efficient,watson2021advancing}); this specific MLE-based solution procedure within the inference-on-PGM family therefore converges to the same classical RSOC policy in the limit. Yet only these \emph{iteration limits} have been identified with interpretable classical objectives; the \emph{surrogate policy} solved for at each iteration is itself a soft, KL-regularized object whose own correspondence to a classical control objective has remained unclear.

Recently, \citet{ito2024risk} addressed this question by replacing KL with R\'enyi divergence in CaI, yielding a family of \emph{log-probability regularized risk-sensitive control} problems with clear, interpretable objectives, indexed by a risk-sensitivity parameter, that in a limiting case recovers the policy obtained by direct conditioning on optimality. Our work offers a complementary, KL-regularization based interpretation, and in doing so addresses two further gaps in the picture. The first is whether the classical SOC and RSOC problems, together with their KL-regularized surrogates from the I- and M-projection formulations, can all be brought together within a single overarching mathematical structure. The second concerns the relation of these KL-regularized formulations to Distributionally Robust Control (DRC), where KL divergences define ambiguity sets for the transition model so as to hedge against model mismatch. Under suitable conditions DRC is known to be equivalent to Risk-Averse Optimal Control with an appropriately chosen risk parameter~(\cite{Nishimura2021_DRO_RSOC, zhangsoft2024, follmer2002convex, noorani2022embracing}), yet how it relates to the other KL-regularized formulations above has remained largely unexplored.

Closely related are the path-integral formulation of \cite{kappen2005linear} and the linearly solvable MDP framework of \cite{todorov2006linearly}, where the control cost itself is a KL divergence between controlled and passive dynamics. The resulting Bellman operator is \emph{linear} and evaluable by sampling the passive system, yielding closed-form policies beyond Linear-Quadratic setting and inspiring sampling-based methods (\cite{williams2019modelTRO}). Generalising these properties is non-trivial; \cite{lefebvre2024probabilistic} showed that the M-projection inherits them, linking these computational properties to the probabilistic graphical model viewpoint.

This paper addresses these gaps by unifying these problems within a broader class of KL-regularized optimal control problems. Prior formulations implicitly impose a \emph{single} KL term acting jointly on policies and transitions, or equivalently, on entire trajectory distributions. We instead introduce a \emph{strict generalization} in which (i) the KL penalties on policies and transitions are separated and (ii) each is weighted independently. This serves as an umbrella formulation, from which classical SOC, RSOC, and their soft-policy variants emerge as special cases. We further show that the soft-policy formulations majorize the original SOC and RSOC objectives via the Majorization--Minimization framework, providing a principled foundation for iterative algorithms. Finally, we identify a \emph{synchronized} sub-case in which the policy and transition KL weights coincide, in which a constellation of remarkable properties emerges simultaneously, revealing that the structural harmony between policy and transition regularization, implicit in many earlier works, is not incidental but mathematically necessary for these properties to hold.

\section{Notation}

We consider a discrete-time controlled stochastic dynamical system over a finite horizon, $T$. The behavior of the system is characterized by trajectories, where a trajectory is the sequence of states and actions, $\underline{\xi}_T = (x_0,u_0,\dots,x_{T-1},u_{T-1},x_T)$.

One way to parametrize a distribution over a system's behavior/trajectories is via a sequence of policies, $\underline{\pi} = (\pi_0,\dots,\pi_{T-1})$, and transition kernels, $\underline{\tau} = (\tau_0,\dots,\tau_{T-1})$. The resulting trajectory distribution then is given by
\begin{equation}\label{eq:traj_dist}
p_{(\underline{\pi},\underline{\tau})}(\underline{\xi}_T) = p(x_0) \prod\nolimits_{t=0}^{T-1} \pi_t(u_t\mid x_t) \,\tau_{t}(x_{t+1}\mid \xi_t),
\end{equation}
where $p(x_0)$ is initial state distribution and $\xi_t=(x_t,u_t)$.

For regularization purposes, we will also consider a baseline behavior, parametrized by a reference policy, $\underline{\rho}$, and the system's transition kernels, $\underline{\iota}$. The corresponding baseline trajectory distribution then is given as
\begin{equation}\label{eq:baseline_traj_dist}
p_{(\underline{\rho},\underline{\iota})}(\underline{\xi}_T) = p(x_0) \prod\nolimits_{t=0}^{T-1} \rho_t(u_t\mid x_t) \,\iota_{t}(x_{t+1}\mid \xi_t).
\end{equation}
Throughout, $\underline{\iota}$ refers to the \emph{system's true transition kernels}---they are not design variables but a fixed reference against which we will measure the deviation of any artificial transition kernel $\underline{\tau}$. The reference policy $\underline{\rho}$ similarly encodes prior or expert behaviour (e.g.\ a stabilising controller or a uniform prior).

Along each trajectory, the system accumulates cost over time. The cumulative cost of a trajectory, $\underline{\xi}_T$, is defined as
\begin{equation}\label{eq:cumulative_cost}
\underline{c}_T(\underline{\xi}_T) = \sum\nolimits_{t=0}^{T-1} c_t(x_t,u_t) + c_T(x_T),
\end{equation}
where $c_t(x_t,u_t) \ge 0$ denotes the running cost incurred by taking action $u_t$ in state $x_t$ at time $t<T$, and $c_T(x_T) \ge 0$ is the terminal cost which depends only on the final state. Throughout, we assume $\underline{c}_T$ is such that all trajectory expectations and risk-sensitive operators arising in the sequel are finite for the policies and transitions considered.

We use the entropic-risk operator $\mathcal{R}^\lambda_{\mu}[f] := -\tfrac{1}{\lambda}\log\mathbb{E}_{\mu}[e^{-\lambda f}]$ and write the KL divergence between $\mu$ and $\sigma$ compactly as $\mathbb{D}^{\mu}_{\sigma}:=\mathbb{E}_{\mu}\!\left[\log\tfrac{\mu}{\sigma}\right]$. Throughout, function arguments are omitted when clear from context: $c_t$, $\pi_t$, $\tau_t$, and the value/$Q$-functions introduced later are to be read with their appropriate domain variables. For a single conditional kernel, $\mathbb{D}^{\pi_t}_{\rho_t}(x_t):=\mathbb{E}_{\pi_t(\cdot|x_t)}[\log\tfrac{\pi_t}{\rho_t}]$ is the standard conditional KL. For sequences $\underline{\pi},\underline{\rho}$ (likewise $\underline{\tau},\underline{\iota}$), with a slight abuse of notation $\mathbb{D}^{\underline{\pi}}_{\underline{\rho}}:=\sum_{t=0}^{T-1}\log\tfrac{\pi_t}{\rho_t}$ denotes the per-trajectory log-ratio sum, always appearing under an outer expectation where, by the chain rule, $\mathbb{E}_{p_{(\underline{\pi},\underline{m})}}[\mathbb{D}^{\underline{\pi}}_{\underline{\rho}}]=\mathbb{D}^{\,p_{(\underline{\pi},\underline{m})}}_{\,p_{(\underline{\rho},\underline{m})}}$ for any common transition kernels $\underline{m}$.

\section{Preliminaries}\label{sec:preliminaries}

In this section we recall two mathematical concepts that will serve as a foundation for developments later in the paper.  
First is the Risk Measures and their dual representations, which we repeatedly use to translate between the Entropic Risk Measure and a KL-regularized expectation-based form.
% and to express certain KL-regularized minimizations in closed form.  
The second is the Majorization--Minimization (MM) framework, which provides the conceptual foundation for later sections, where KL-regularized problems are used as tractable surrogates for classical control objectives.

\subsection{Risk Measures and their Dual Representations} \label{subsec:risk_measures_and_duals}

Risk measures, originating in mathematical finance \cite{follmer2002convex}, quantify preferences over uncertain outcomes beyond a plain expectation and thereby encode risk-sensitive behaviour. Through their \emph{dual representation}, every (convex) risk measure can be written as a worst-/best-case expectation over a set of alternative distributions penalised by a divergence from a reference distribution $\rho$~\cite{follmer2002convex}. This duality is the bridge between risk sensitivity and KL regularisation that we exploit throughout.

Among the many available risk measures, in this paper we work exclusively with the \emph{entropic risk measure}, whose KL-based dual representation gives, for $\lambda\in\mathbb{R}\setminus\{0\}$,
\begin{equation}\label{eq:generalized-risk-operator}
\mathcal{R}^{\lambda}_{\rho}[f] = 
\begin{cases}
\displaystyle\inf_{\pi \ll \rho} \left\{ \mathbb{E}_{\pi}[f] + \tfrac{1}{\lambda} \mathbb{D}^{\pi}_{\rho} \right\}, & \lambda > 0\;(\text{risk-seeking}),\\[2pt]
\displaystyle\sup_{\pi \ll \rho} \left\{ \mathbb{E}_{\pi}[f] + \tfrac{1}{\lambda} \mathbb{D}^{\pi}_{\rho} \right\}, & \lambda < 0\;(\text{risk-averse}),
\end{cases}
\end{equation}
with the same closed-form extremiser in both cases,
\begin{equation}\label{eq:generalized-optimal-dist}
\pi^* = \frac{\rho\, e^{-\lambda f}}{\mathbb{E}_{\rho}[e^{-\lambda f}]}.
\end{equation}

The dual representations \eqref{eq:generalized-risk-operator} provide a crucial link between risk sensitivity and regularization: the risk assessment of a cost function considers alternative distributions $\pi$, where deviations from $\rho$ are penalized by the KL divergence, along with the direction of optimization (minimization vs. maximization) determining the risk attitude. This perspective also establishes a connection to DRC with KL ambiguity sets. Here a saddle point problem is solved to calculate the exact Lagrangian multiplier that corresponds with a predefined KL discrepancy. In that sense, Risk Sensitivity can be interpreted as a soft DRC problem where the Lagrangian multiplier is predefined instead.

To capture both risk attitudes in a single notation, we use the sign-dependent extremiser
\[
\operatorname*{opt}^{\lambda}_{q \in \mathcal{Q}} J(q) := 
\begin{cases}
\displaystyle \min_{q \in \mathcal{Q}} J(q), & \lambda > 0,\\
\displaystyle \max_{q \in \mathcal{Q}} J(q), & \lambda < 0,
\end{cases}
\]
with $\arg\operatorname*{opt}^{\lambda}$ the corresponding extremiser.

\subsection{Majorization--Minimization (MM) Framework}

The Majorization--Minimization framework addresses optimization problems of the form $\min_{x} F(x)$ particularly when $F$ is difficult to optimize directly. MM algorithms iteratively minimize tractable surrogate functions that upper-bound $F$.

Given an iterate $x^{k}$, a function $G(\cdot\,|\,x^{k})$ \emph{majorizes} $F$ at $x^{k}$ if
\begin{subequations}
    \begin{align}
        G(x^{k}\,|\,x^{k}) &= F(x^{k}), \\
        G(x\,|\,x^{k}) &\ge F(x) \;\forall x,
    \end{align}
\end{subequations}
and the next iterate is defined by $x^{k+1} = \arg\min_x G(x\,|\,x^{k})$.

Two features make MM attractive:
\begin{enumerate}
\item \textbf{Descent guarantee:}
Since minimization of $G$ yields $G(x^{k+1}|x^{k})\le G(x^{k}|x^{k})$, majorization implies that $F(x^{k+1}) \le F(x^{k}).$
\item \textbf{Surrogate flexibility:}
The surrogate $G$ may be chosen to be a simpler tractable objective, allowing efficient updates.
\end{enumerate}

In later sections, we will see that certain KL-regularized control problems naturally provide such surrogates for classical optimal control problems. This facilitates algorithms in which each iteration solves a tractable KL-regularized subproblem while guaranteeing descent on the original objective.

\section{A Unifying KL-Regularized Control Formulation} \label{sec:unifying_kl}

We introduce a central KL-regularized optimal control problem (C-KLR-OC) that jointly optimises a policy sequence $\underline{\pi}$ and an \emph{artificial} transition-kernel sequence $\underline{\tau}$, while penalising the KL deviation of each from baseline behaviour $(\underline{\rho},\underline{\iota})$. Promoting the transition kernels to decision variables endows the design with the capacity to reason about alternative system evolutions, whether optimistic or pessimistic, thereby encoding risk-sensitive behaviour.

The central KL-regularized optimal control problem is defined as:
\begin{equation}\label{eq:central_problem_full}
\begin{aligned}
\min_{\underline{\pi}}
\operatorname*{opt}^{\lambda^S}_{\underline{\tau}}
\mathbb{E}_{p_{(\underline{\pi},\underline{\tau})}}\!\Bigl[\underline{c}_T
+ \frac{1}{\lambda^{P}}\mathbb{D}^{\underline{\pi}}_{\underline{\rho}}
+ \frac{1}{\lambda^{S}}\mathbb{D}^{\underline{\tau}}_{\underline{\iota}}
\Bigr],
\end{aligned}
\end{equation}
where $\lambda^{P}>0$ controls deviations of the policy from $\underline{\rho}$ and $\lambda^{S}\in\mathbb{R}\setminus\{0\}$ controls deviations of the transition kernels from $\underline{\iota}$, with the sign of $\lambda^{S}$ selecting an optimistic ($\lambda^S>0$, risk-seeking) or pessimistic ($\lambda^S<0$, risk-averse) attitude. We focus on the case in which $\lambda^P$ and $\lambda^S$ are constant across stages; an extension to non-constant weights is left to future work.

\begin{thm}[Optimal Solution for C-KLR-OC] 
\label{thm:central_kl_problem} 
For the problem defined in \eqref{eq:central_problem_full}, define the terminal value function $V_T = c_T$. Then, for $t = T-1, \dots, 0$, the optimal value function $V_t$, action-value function $Q_t$, policy $\pi_t^*$, and transition kernel $\tau_t^*$ satisfy the following recursions:

\begin{subequations}\label{eq:central_recursion}
\begin{align}
V_t &= \min_{\pi_t} \mathbb{E}_{\pi_t}[Q_t] + \frac{1}{\lambda^P} \mathbb{D}^{\pi_t}_{\rho_t} = \mathcal{R}^{\lambda^P}_{\rho_t}[Q_t], \\
Q_t &= \operatorname*{opt}^{\lambda^S}_{\tau_{t}} \left\{ \mathbb{E}_{\tau_{t}}[c_t + V_{t+1}] + \frac{1}{\lambda^S} \mathbb{D}^{\tau_t}_{\iota_t} \right\} = \mathcal{R}^{\lambda^S}_{\iota_t}[c_t+V_{t+1}], \\
\pi^*_t &= \arg\min_{\pi_t} \mathbb{E}_{\pi_t}[Q_t] + \frac{1}{\lambda^P} \mathbb{D}^{\pi_t}_{\rho_t} = \rho_t \frac{e^{-\lambda^P Q_t}}{e^{-\lambda^P V_t}}, \\
\tau^*_{t} &= \arg\operatorname*{opt}^{\lambda^S}_{\tau_{t}} \left\{ \mathbb{E}_{\tau_{t}}[V_{t+1}] + \frac{1}{\lambda^S} \mathbb{D}^{\tau_t}_{\iota_t} \right\} = \iota_t \frac{e^{-\lambda^S (c_t+V_{t+1})}}{e^{-\lambda^S Q_t}}.
\end{align}
\end{subequations}
\end{thm}

\begin{pf}
See Appendix~\ref{proof:central_kl_problem}.
\end{pf}

\begin{rem}
This formulation represents a departure from the classical Probabilistic Control Design (PCD) paradigm \cite{karny1996towards}, where one begins with a prescribed optimal trajectory distribution and seeks policies that \emph{match it}. In our approach, we do not presuppose the form of the optimal trajectory distribution; instead, it emerges organically as the solution to the optimization problem \eqref{eq:central_problem_full}.     
\end{rem}

\begin{rem}
The computational tractability of the central KL-regularized problem \eqref{eq:central_problem_full} stems from a fundamental shift in the optimization geometry over probability measures. In the unregularized case, the inner optimizations $\min_{\pi_t} \mathbb{E}_{\pi_t}[Q_t]$ and $\operatorname*{opt}_{\tau_t} \mathbb{E}_{\tau_t}[V_{t+1}]$ over unrestricted measures invariably yield Diracs concentrated at the global extrema $\arg\min_u$ and $\arg\max_{x'}$ of the underlying continuous spaces, requiring computationally prohibitive, generally non-convex global searches. The KL penalties regularize these measure-space optimizations, replacing the point-mass solutions with closed-form Gibbs measures (Theorem~\ref{thm:central_kl_problem}). Substituting these back into the Bellman recursion analytically converts the non-differentiable $\min$ and $\max$ operators into the smooth entropic risk operator, so the intractable search for global extrema becomes the evaluation of smooth expectations, which readily admits scalable numerical approximations including the Majorization-Minimization framework of Section~\ref{sec:MM}.
\end{rem}

\section{Recovering Standard Control Formulations}
\label{sec:recovering_formulations}

The central problem~\eqref{eq:central_problem_full} acts as an umbrella under which classical formulations are recovered by toggling (i)~whether the policy is regularised against a baseline ($\underline{\rho}\neq\underline{\pi}$) or not ($\underline{\rho}=\underline{\pi}$), and (ii)~whether the transition kernels are fixed to $\underline{\iota}$ or left free. Table~\ref{tab:dp_recursions} summarises the resulting dynamic programming recursions.

\subsection{No Policy Regularization: $\underline{\rho} = \underline{\pi}$}

When no meaningful baseline policy is available (e.g., when designing a policy from scratch), the policy regularization term can be removed by setting the reference policy equal to the optimizing policy, yielding $\mathbb{D}^{\pi_t}_{\rho_t}=0$. Two cases arise depending on whether the transition kernels are fixed/free.

\subsubsection{Stochastic Optimal Control (SOC)}\label{sec:SOC}
Setting $\underline{\rho}=\underline{\pi}$ and constraining $\underline{\tau}=\underline{\iota}$, the central viewpoint yields the good old SOC problem
\begin{equation}
\label{eq:SOC_min}
\min_{\underline{\pi}} \; \mathbb{E}_{p_{(\underline{\pi},\underline{\iota})}}[\underline{c}_T].
\end{equation}

\subsubsection{Risk-Sensitive Stochastic Optimal Control (RSOC)}\label{sec:RSOC}
By setting $\underline{\rho} = \underline{\pi}$, while leaving $\underline{\tau}$ free, the central problem becomes
\begin{equation}
\label{eq:RSOC_full}
\min_{\underline{\pi}} \operatorname*{opt}^{\lambda^S}_{\underline{\tau}}
\mathbb{E}_{p_{(\underline{\pi},\underline{\tau})}}[\underline{c}_T]
+ \tfrac{1}{\lambda^S}
\mathbb{D}^{p_{(\underline{\pi},\underline{\tau})}}_{p_{(\underline{\pi},\underline{\iota})}}.
\end{equation}
Equation \eqref{eq:RSOC_full} admits an equivalent two-step representation:
\begin{subequations}
  \begin{align}
&\min_{\underline{\pi}} \operatorname*{opt}^{\lambda^S}_{\underline{\tau}}
\mathbb{E}_{p_{(\underline{\pi},\underline{\tau})}}[\underline{c}_T]
+ \tfrac{1}{\lambda^S}
\mathbb{D}^{p_{(\underline{\pi},\underline{\tau})}}_{p_{(\underline{\pi},\underline{\iota})}} \notag
\\
&=
\min_{\underline{\pi}}
\Bigg\{
\operatorname*{opt}^{\lambda^S}_{p(\underline{\xi}_T)}\;
\mathbb{E}_{p(\underline{\xi}_T)}[\underline{c}_T]
+
\tfrac{1}{\lambda^S}\mathbb{D}^{p(\underline{\xi}_T)}_{p_{(\underline{\pi},\underline{\iota})}}
\Bigg\} \label{eq:RS_rep2} \\
&=
\min_{\underline{\pi}}
\left[
-\tfrac{1}{\lambda^S}\log
\mathbb{E}_{p_{(\underline{\pi},\underline{\iota})}}
\left[e^{-\lambda^S\underline{c}_T}\right]
\right]. \label{eq:RS_rep3}
\end{align}  
\end{subequations}
Equation \eqref{eq:RS_rep3} is the standard risk-sensitive stochastic optimal control formulation. When $\lambda^S > 0$, it yields risk-seeking behavior, where the agent exhibits optimistic behavior by optimizing against favorable transitions. When $\lambda^S < 0$, it yields risk-averse behavior, where the agent exhibits pessimistic behavior by optimizing against worst-case transitions.

\subsection{With Policy Regularization: $\underline{\rho} \neq \underline{\pi}$}

When a meaningful baseline $\underline{\rho}$ is available (expert demonstrations, safe prior, stabilising controller), regularising towards it is informative. Fixing or freeing $\underline{\tau}$ yields two formulations, whose optima are \emph{soft} policies, that is, exponentially tilted, cost-favoured versions of $\underline{\rho}$.

\subsubsection{Soft-Policy SOC (SP-SOC)}
Fixing $\underline{\tau}=\underline{\iota}$ and taking $\underline{\rho}$ as a baseline policy, the central problem becomes
\begin{equation}
\label{eq:KLR_SOC}
\min_{\underline{\pi}}
\mathbb{E}_{p_{(\underline{\pi},\underline{\iota})}}[\underline{c}_T]
+ \frac{1}{\lambda^P}
\mathbb{D}^{p_{(\underline{\pi},\underline{\iota})}}_{p_{(\underline{\rho},\underline{\iota})}}.
\end{equation}

\begin{rem}
This formulation can be interpreted as a density-matching problem. Specifically, it is equivalent to the I-projection problem \cite{neumann2011variational,lefebvre2024probabilistic}:
\begin{equation}
\min_{\underline{\pi}} \mathbb{D}^{p_{(\underline{\pi},\underline{\iota})}}_{p^*}
\end{equation}
where $p^* \propto p_{(\underline{\rho},\underline{\iota})} e^{-\lambda^P \underline{c}_T}$ is the target distribution. 

Note that when $\underline{\rho}$ is uniform, the KL term reduces to policy entropy, recovering the maximum entropy SOC formulation.
\end{rem}

\subsubsection{Soft-Policy Risk-Sensitive Stochastic Optimal \\Control (SP-RSOC)}
Allowing both $\underline{\pi}$ and $\underline{\tau}$ to be optimized with regularization parameters $\lambda^P > 0$ and $\lambda^S \neq 0$ produces:
\begin{equation}
\label{eq:KLR_RSOC_general}
\min_{\underline{\pi}} \operatorname*{opt}^{\lambda^S}_{\underline{\tau}}
\mathbb{E}_{p_{(\underline{\pi},\underline{\tau})}}\!\Bigl[\underline{c}_T
+ \frac{1}{\lambda^P}
\mathbb{D}^{\underline{\pi}}_{\underline{\rho}}
+ \frac{1}{\lambda^S}
\mathbb{D}^{\underline{\tau}}_{\underline{\iota}}\Bigr].
\end{equation}
Here $\lambda^S > 0$ a risk-seeking behavior, while $\lambda^S < 0$ induces risk-averse behavior. 

When $\lambda^P = \lambda^S = \lambda > 0$, we get:
\begin{equation}
\label{eq:KLR_RSOC_sync}
\min_{\underline{\pi}} \min_{\underline{\tau}}
\mathbb{E}_{p_{(\underline{\pi},\underline{\tau})}}\!\Bigl[\underline{c}_T
+ \frac{1}{\lambda}
\mathbb{D}^{\underline{\pi}}_{\underline{\rho}}
+ \frac{1}{\lambda}
\mathbb{D}^{\underline{\tau}}_{\underline{\iota}}\Bigr].
\end{equation}
We denote this case as the \emph{Synchronized SP-RSOC} (S-SP-RSOC). Note that as $\lambda^S>0$, this formulation has a risk-seeking attitude. It exhibits several remarkable properties which we detail in Section~\ref{sec:sp_rsoc_features}.

\begin{rem}
The S-SP-RSOC problem \eqref{eq:KLR_RSOC_sync} admits the same solution as a density-matching problem based on the \emph{M-projection} \cite{neumann2011variational,lefebvre2024probabilistic}:
\begin{equation}
    \min_{\underline{\pi}} \mathbb{D}^{p^*}_{p_{(\underline{\pi},\underline{\iota})}},
\end{equation}
where $p^* \propto p_{(\underline{\rho},\underline{\iota})} e^{-\lambda \underline{c}_T}$, is the target distribution.
\end{rem}

\begin{rem}
By replacing KL with R\'enyi divergence in CaI variational inference, \citet{ito2024risk} derived a family of \emph{log-probability regularized risk-sensitive control} problems, indexed by a risk-sensitivity parameter, whose optima in a limiting case recover the policy obtained by direct conditioning on optimality. Our S-SP-RSOC~\eqref{eq:KLR_RSOC_sync} offers a complementary, KL-regularization based interpretation of the same direct-conditioning policy as the solution of a risk-sensitive problem in which policy \emph{and} transitions are joint design variables, each KL-penalised against its baseline.
\end{rem}

\subsection{Equivalences for Deterministic Dynamics}
\label{subsec:deterministic_dynamics}

When the baseline dynamics $\underline{\iota}$ are \emph{deterministic} (each $\iota_t$ a Dirac at the prescribed next state), absolute continuity $\tau_t \ll \iota_t$ forces every feasible $\tau_t$ to coincide with $\iota_t$, so the auxiliary transition kernels of RSOC and SP-RSOC lose all expressive freedom and the $\underline{\tau}$-extremisation reduces to the same constraint that defines SOC and SP-SOC. Consequently,
\[
\text{SOC} \equiv \text{RSOC},\qquad \text{SP-SOC} \equiv \text{SP-RSOC},
\]
and we refer to these deterministic special cases as Deterministic Optimal Control (DOC) and Soft-Policy DOC (SP-DOC).

\section{Majorization of Standard Control Formulations}\label{sec:MM}

We now show that SP-SOC and SP-RSOC serve as majorizers for their classical counterparts, SOC and RSOC respectively.   
This provides a principled interpretation: KL-regularized problems are not merely relaxations, but tractable surrogate objectives guaranteeing descent on the corresponding control objectives.

\subsection{SP-SOC majorizes SOC.}

For SOC, we have the objective: $J_{\mathrm{SOC}}(\underline{\pi}) = \mathbb{E}_{p_{(\underline{\pi},\underline{\iota})}}[\underline{c}_T]$. Given an iterate $\underline{\pi}^{k}$, the SP-SOC surrogate with $\lambda^P>0$ and baseline $\underline{\rho}=\underline{\pi}^{k}$ reduces to
\begin{equation}
J_{\mathrm{SP-SOC}}(\underline{\pi}|\underline{\pi}^{k})
:=
\mathbb{E}_{p_{(\underline{\pi},\underline{\iota})}}\!\Bigl[\underline{c}_T
+ \tfrac{1}{\lambda^P}
\mathbb{D}^{\underline{\pi}}_{\underline{\pi}^{k}}\Bigr].
\end{equation}
By the chain rule of KL, $J_{\mathrm{SP-SOC}}(\underline{\pi}|\underline{\pi}^{k})-J_{\mathrm{SOC}}(\underline{\pi}) = \tfrac{1}{\lambda^P}\mathbb{E}_{p_{(\underline{\pi},\underline{\iota})}}[\mathbb{D}^{\underline{\pi}}_{\underline{\pi}^{k}}] = \tfrac{1}{\lambda^P}\mathbb{D}^{\,p_{(\underline{\pi},\underline{\iota})}}_{\,p_{(\underline{\pi}^{k},\underline{\iota})}} \ge 0$, with equality iff $\underline{\pi}=\underline{\pi}^{k}$. This establishes the MM conditions $J_{\mathrm{SP-SOC}}(\underline{\pi}^{k}|\underline{\pi}^{k})=J_{\mathrm{SOC}}(\underline{\pi}^{k})$ and $J_{\mathrm{SP-SOC}}(\underline{\pi}|\underline{\pi}^{k})\ge J_{\mathrm{SOC}}(\underline{\pi})$, so $\underline{\pi}^{k+1}\leftarrow \arg\min_{\underline{\pi}} J_{\mathrm{SP-SOC}}(\underline{\pi}|\underline{\pi}^{k})$ guarantees descent on $J_{\mathrm{SOC}}$.

\subsection{SP-RSOC majorizes RSOC.}

Starting from RSOC~\eqref{eq:RSOC_full} and cancelling $\underline{\pi}$ inside the trajectory KL gives the equivalent form $J_{\mathrm{RSOC}}(\underline{\pi},\underline{\tau})=\mathbb{E}_{p_{(\underline{\pi},\underline{\tau})}}\!\bigl[\underline{c}_T+\tfrac{1}{\lambda^S}\mathbb{D}^{\underline{\tau}}_{\underline{\iota}}\bigr]$. With $\underline{\pi}^{k}$ as baseline and $\lambda^P>0$, the SP-RSOC surrogate
\begin{equation}
\label{eq:J_KLR_RSOC}
J_{\mathrm{SP-RSOC}}(\underline{\pi},\underline{\tau}|\underline{\pi}^{k})
:=
\mathbb{E}_{p_{(\underline{\pi},\underline{\tau})}}\!\Bigl[\underline{c}_T
+ \tfrac{1}{\lambda^P}
\mathbb{D}^{\underline{\pi}}_{\underline{\pi}^{k}}
+ \tfrac{1}{\lambda^S}
\mathbb{D}^{\underline{\tau}}_{\underline{\iota}}\Bigr]
\end{equation}
differs from $J_{\mathrm{RSOC}}$ by the policy-KL term alone. By the chain rule of KL, $J_{\mathrm{SP-RSOC}}(\underline{\pi},\underline{\tau}|\underline{\pi}^{k})-J_{\mathrm{RSOC}}(\underline{\pi},\underline{\tau}) = \tfrac{1}{\lambda^P}\mathbb{D}^{\,p_{(\underline{\pi},\underline{\tau})}}_{\,p_{(\underline{\pi}^{k},\underline{\tau})}} \ge 0$ for every $(\underline{\pi},\underline{\tau})$, with equality iff $\underline{\pi}=\underline{\pi}^{k}$. As in the previous subsection, the MM conditions follow for the outer $\underline{\pi}$-objective, so the iteration $\underline{\pi}^{k+1}\leftarrow\arg\min_{\underline{\pi}}\operatorname*{opt}^{\lambda^S}_{\underline{\tau}} J_{\mathrm{SP-RSOC}}(\underline{\pi},\underline{\tau}|\underline{\pi}^{k})$ guarantees descent on $J_{\mathrm{RSOC}}$.

\section{Special Properties of S-SP-RSOC}
\label{sec:sp_rsoc_features}

The S-SP-RSOC formulation exhibits several remarkable structural properties that distinguish it from the other classical formulations and enable powerful extensions such as compositional control design. Note that they only hold when the policy and transition regularization weights are equal ($\lambda^P = \lambda^S = \lambda > 0$). These properties have previously been identified in continuous-time Path Integral Control~\cite{kappen2005linear} and in Linear Markov Games and the SP-DOC case~\cite{dvijotham2012linearly}; our analysis reveals that they extend to the broader class of S-SP-RSOC.

\subsection{Linear Bellman Operator}

For S-SP-RSOC, \eqref{eq:central_recursion} reduces to
\begin{subequations}
\begin{align}
V_t &= \mathcal{R}^{\lambda}_{\rho_t}[Q_t] = -\tfrac{1}{\lambda}\log\mathbb{E}_{\rho_t}\!\left[e^{-\lambda Q_t}\right], \\
Q_t &= \mathcal{R}^{\lambda}_{\iota_t}[c_t+V_{t+1}] = -\tfrac{1}{\lambda}\log\mathbb{E}_{\iota_t}\!\left[e^{-\lambda(c_t+V_{t+1})}\right].
\end{align}
\end{subequations}
Introducing the \emph{desirability} $z_t:=e^{-\lambda V_t}$ and \emph{reward} $r_t:=e^{-\lambda c_t}$, the Bellman equation becomes linear:
\begin{equation}\label{eq:sp_rsoc_linear}
z_t = \mathbb{E}_{\rho_t}\left[r_t\,\mathbb{E}_{\iota_t}\!\left[z_{t+1}\right]\right].
\end{equation}
This linearity arises from the multiplicative structure of the exponential transform and represents a significant simplification compared to the nonlinear Bellman equations of standard stochastic optimal control.

\subsection{Path Integral Solution}

The linearity of \eqref{eq:sp_rsoc_linear} unrolls into a \emph{path integral}. Specifically, the desirability is an expectation over trajectories generated by the baseline,
\begin{equation}\label{eq:sp_rsoc_path_integral}
z_t = \mathbb{E}_{p_{(\underline{\rho}_t,\underline{\iota}_t|x_t)}}\!\left[e^{-\lambda\left(\sum_{k=t}^{T-1} c_k + c_T\right)}\right],
\end{equation}
where $p_{(\underline{\rho}_t,\underline{\iota}_t|x_t)}$ is the trajectory distribution under the baseline $(\underline{\rho}_t,\underline{\iota}_t)$ from time $t$ onward, conditioned on $x_t$. The value function can therefore be estimated by \emph{forward} Monte-Carlo sampling under the baseline. This approach bypasses backward DP, is naturally parallel, and is model-free given the baseline. The optimal policy is the baseline reweighted by reward and expected future desirability,
\begin{equation}
\pi_t^* = \rho_t\,\frac{r_t\,\mathbb{E}_{\iota_t}[z_{t+1}]}{z_t}.
\end{equation}

\subsection{Compositionality of Value Functions and Policies}

The linearity of the desirability recursion immediately yields compositionality, in the sense that solutions of simpler subproblems combine into solutions of more complex ones.

\begin{thm}[Compositionality of S-SP-RSOC]\label{th:CompositionalityOfSRS-SP-RSOC}
Suppose the terminal cost decomposes as
\begin{equation}\label{eq:compositional_terminal}
e^{-\lambda c_T} = \sum\nolimits_{n=1}^N \gamma_n e^{-\lambda c_T^{(n)}},
\end{equation}
with weights $\gamma_n > 0$, and for each $n = 1,\dots,N$ define the component desirabilities recursively by
\begin{subequations}\label{eq:component_desirability}
\begin{align}
z_T^{(n)} &:= \gamma_n e^{-\lambda c_T^{(n)}}, \\
z_t^{(n)} &:= \mathbb{E}_{\rho_t}\!\left[r_t\,\mathbb{E}_{\iota_t}\!\left[z_{t+1}^{(n)}\right]\right], \quad t = T-1,\dots,0.
\end{align}
\end{subequations}
Then, for all $t = 0,\dots,T-1$, the desirability and the optimal policy admit the decompositions
\begin{equation}
z_t = \sum\nolimits_{n=1}^N z_t^{(n)},
\qquad
\pi_t^* = \sum\nolimits_{n=1}^N w_t^{(n)} \pi_t^{(n)},
\end{equation}
with mixture weights $w_t^{(n)} = z_t^{(n)}/z_t$, $\sum_{n=1}^N w_t^{(n)} = 1$, and component policies
\begin{equation}
\pi_t^{(n)} = \rho_t\,\frac{r_t\,\mathbb{E}_{\iota_t}[z_{t+1}^{(n)}]}{z_t^{(n)}}.
\end{equation}
\end{thm}

\begin{pf}
See Appendix~\ref{proof:compsitionality}.
\end{pf}

\subsection{Maximum Likelihood Estimation on a PGM}

We conclude by establishing a final connection with the CaI framework \cite{toussaint2006probabilistic,levine2018reinforcement,lefebvre2024probabilistic}, which introduces binary optimality variables $\mathcal{O}_t\in\{0,1\}$ on the system PGM with likelihood $p(\underline{\mathcal{O}}_T=\underline{1}_T|\underline{\xi}_T)\propto e^{-\lambda\underline{c}_T}$, $\lambda>0$. Several connections to our framework follow.

\textbf{Posterior equivalence.} For S-SP-RSOC with $\lambda^P=\lambda^S=\lambda>0$, denoting $q:=p_{(\underline{\pi},\underline{\tau})}$ and $p_b:=p_{(\underline{\rho},\underline{\iota})}$, the central problem reduces to $\min_{q\ll p_b}\{\mathbb{E}_q[\underline{c}_T]+\tfrac{1}{\lambda}\mathbb{D}^{q}_{p_b}\}$. By~\eqref{eq:generalized-optimal-dist},
\begin{equation}
q^* = \frac{p_b\, e^{-\lambda\underline{c}_T}}{\mathbb{E}_{p_b}[e^{-\lambda\underline{c}_T}]} = p(\underline{\xi}_T\mid\underline{\mathcal{O}}_T=1;\underline{\rho},\underline{\iota}),
\end{equation}
i.e. the CaI posterior of states and actions given optimality.

\textbf{Density matching.} Using $q^*$ as a target, forward-KL minimisation $\min_{\underline{\pi}}\mathbb{D}^{p_{(\underline{\pi},\underline{\iota})}}_{q^*}$ recovers SP-SOC; reverse-KL minimisation $\min_{\underline{\pi}}\mathbb{D}^{q^*}_{p_{(\underline{\pi},\underline{\iota})}}$ recovers S-SP-RSOC.

\textbf{MLE interpretation.} Maximising the likelihood of observing optimality: $\max_{\underline{\pi}}\,\log p(\underline{\mathcal{O}}_T=\underline{1}_T;\underline{\pi},\underline{\iota})$, can be shown to be equivalent to the risk-seeking RSOC problem~\eqref{eq:RS_rep3} and can be solved efficiently by the Expectation--Maximisation algorithm, whose E-step computes the posterior $q^k(\underline{\xi}_T)=p(\underline{\xi}_T\mid\underline{\mathcal{O}}_T;\underline{\pi}^k)$ and whose M-step updates the policy as $\underline{\pi}^{k+1} = \arg\max_{\underline{\pi}}\,\mathbb{E}_{q^k}[\log p(\underline{\xi}_T,\underline{\mathcal{O}}_T;\underline{\pi})]$.

\textbf{Conditional policy form.} The optimal S-SP-RSOC policy is exactly $\pi_t^\text{S-SP-RSOC}(u_t|x_t)=p(u_t|x_t,\underline{\mathcal{O}}_T=1;\underline{\rho},\underline{\iota})$, computable by Bayesian smoothing, which is why the S-SP-RSOC recursions take the form of backward message passing. This identification with exact PGM inference explains the linear Bellman, path-integral, and compositionality properties.
% \vspace{-6pt}

\section{Conclusion}\label{sec:conclusion}

This paper consolidated several optimal-control formulations through a single Central KL-regularized Optimal Control (C-KLR-OC) problem in which the policy and transition KL penalties are separated and weighted independently. From this umbrella formulation, SOC, RSOC, SP-SOC, and SP-RSOC arise as structural restrictions, and the soft-policy formulations act as surrogates yielding descent-guaranteeing iterations on the classical objectives. A key insight is that the policy and transition regularisations need not be equivalent. Put differently, the convergence rate of the soft-policy fixed-point iteration and the risk parameter of the corresponding RSOC need not be synchronised. When they are, the resulting S-SP-RSOC exhibits favourable computational properties including linear Bellman operator, path-integral solution, and compositionality, showing this harmony to be mathematically necessary, not incidental, for these properties to hold simultaneously beyond previously known settings.

Our future work will explore non-constant regularization weights $\lambda^P,\lambda^S$ for finer control; replacing soft regularization by hard KL constraints; integrating partial observability; and translating these structural insights into practical algorithms.

\bibliography{ifacconf}             % bib file to produce the bibliography

\appendix
\appendix

\section{Proof of Theorem \label{proof:central_kl_problem}}
\label{app:proof_central_problem}

We prove the dynamic programming recursion for the central KL-regularized optimal control problem \eqref{eq:central_problem_full} by backward induction on the time index $t$.

\textbf{Notation:}
For $t = 0, 1, \dots, T-1$, let:
\begin{itemize}
    \item $\underline{\pi}_t = (\pi_t, \pi_{t+1}, \dots, \pi_{T-1})$ denote the policy sequence from time $t$ onward
    \item $\underline{\tau}_t = (\tau_t, \tau_{t+1}, \dots, \tau_{T-1})$ denote the transition kernel sequence from time $t$ onward
    \item $p_{(\underline{\pi}_t,\underline{\tau}_t)}$ denote the conditional distribution of the trajectory segment from $t$ to $T$, given $x_t$
\end{itemize}

Define the cost-to-go from state $x_t$ under policy and transition sequences $(\underline{\pi}_t, \underline{\tau}_t)$ as:

\begin{equation}
\adjustbox{max width=\linewidth}{%
$\displaystyle J_t(x_t; \underline{\pi}_t, \underline{\tau}_t) = \mathbb{E}_{p_{(\underline{\pi}_t,\underline{\tau}_t)}} \left[ c_T + \sum_{k=t}^{T-1} \left( c_{k} + \frac{1}{\lambda^P} \mathbb{D}^{\pi_{k}}_{\rho_{k}} + \frac{1}{\lambda^S} \mathbb{D}^{\tau_{k}}_{\iota_{k}} \right) \right]$%
}
\end{equation}

The optimal cost-to-go function is:
\begin{equation}
V_t(x_t) = \min_{\underline{\pi}_t} \operatorname*{opt}^{\lambda^S}_{\underline{\tau}_t} J_t(x_t; \underline{\pi}_t, \underline{\tau}_t),
\end{equation}

\textbf{Proof Structure:}
We will prove by backward induction that $V_t$ satisfies the recursions in \eqref{eq:central_recursion}. The proof proceeds in three stages:

\begin{enumerate}
    \item \textbf{Base Case ($t=T$)}: Define the terminal condition.
    \item \textbf{Inductive Hypothesis}: Assume the recursion holds at time $t+1$.
    \item \textbf{Inductive Step}: Prove the recursion holds at time $t$.
\end{enumerate}

\textbf{Base Case ($t=T$):}
At the terminal time $T$, there are no decisions to make. By definition:
\begin{equation}
V_T(x_T) = c_T(x_T).
\end{equation}
This serves as the initial condition for the backward recursion. 

\textbf{Inductive Hypothesis:}
For the inductive proof, we assume that for time $t+1$ (where $0 \leq t+1 \leq T-1$), the function $V_{t+1}$ exists and satisfies the recursions \eqref{eq:central_recursion}. Specifically, we assume:

\begin{subequations}
\begin{align}
V_{t+1} &= \min_{\pi_{t+1}} \left[ \mathbb{E}_{\pi_{t+1}}[Q_{t+1}] + \frac{1}{\lambda^P} \mathbb{D}^{\pi_{t+1}}_{\rho_{t+1}} \right] = \mathcal{R}^{\lambda^P}_{\rho_{t+1}}[Q_{t+1}], \\
Q_{t+1} &= \operatorname*{opt}^{\lambda^S}_{\tau_{t+1}} \left\{ \mathbb{E}_{\tau_{t+1}}[c_{t+1} + V_{t+2}] + \frac{1}{\lambda^S} \mathbb{D}^{\tau_{t+1}}_{\iota_{t+1}} \right\} \notag \\ 
&= \mathcal{R}^{\lambda^S}_{\iota_{t+1}}[c_{t+1} + V_{t+2}], \\
\pi_{t+1}^* &= \rho_{t+1} \frac{e^{-\lambda^P Q_{t+1}}}{e^{-\lambda^P V_{t+1}}}, \\
\tau_{t+1}^* &= \iota_{t+1} \frac{e^{-\lambda^S (c_{t+1} + V_{t+2})}}{e^{-\lambda^S Q_{t+1}}}.
\end{align}
\end{subequations}

This hypothesis is trivially true for $t+1 = T-1$ when we take $V_T$ as defined in base case, and use dual representations \eqref{eq:generalized-risk-operator} and the extremal distribution \eqref{eq:generalized-optimal-dist}.

\textbf{Inductive Step:}
We now prove that if the inductive hypothesis holds for $t+1$, then the recursions \eqref{eq:central_recursion} hold for $t$.

By the principle of optimality (Bellman's principle), an optimal policy from time $t$ onward must consist of an optimal immediate decision at time $t$ followed by optimal decisions from time $t+1$ onward. Therefore:

\begin{equation}\label{eq:bellman_opt}
\begin{aligned}
V_t(x_t) &= \min_{\pi_t} \operatorname*{opt}^{\lambda^S}_{\tau_t} \mathbb{E}_{\pi_t, \tau_t} \Bigg[ c_t + \frac{1}{\lambda^P} \mathbb{D}^{\pi_{t}}_{\rho_{t}} \\
&\quad + \frac{1}{\lambda^S} \mathbb{D}^{\tau_{t}}_{\iota_{t}} + \underbrace{\min_{\underline{\pi}_{t+1}} \operatorname*{opt}^{\lambda^S}_{\underline{\tau}_{t+1}} J_{t+1}(x_{t+1}; \underline{\pi}_{t+1}, \underline{\tau}_{t+1})}_{= V_{t+1}(x_{t+1})} \Bigg].
\end{aligned}
\end{equation}

For fixed $(x_t, u_t)$, define $f(x_{t+1}) = c_t(x_t, u_t) + V_{t+1}(x_{t+1})$. The term involving $\tau_t$ in \eqref{eq:bellman_opt} is:

\begin{equation}
\operatorname*{opt}^{\lambda^S}_{\tau_t} \left\{ \mathbb{E}_{\tau_t}[f] + \frac{1}{\lambda^S} \mathbb{D}^{\tau_{t}}_{\iota_{t}} \right\}.
\end{equation}

This is exactly of the form in the dual representation of the entropic risk measure \eqref{eq:generalized-risk-operator}. Therefore, using this dual representation:

\begin{equation}
\operatorname*{opt}^{\lambda^S}_{\tau_t} \left\{ \mathbb{E}_{\tau_t}[f] + \frac{1}{\lambda^S} \mathbb{D}^{\tau_t}_{\iota_t} \right\} = \mathcal{R}^{\lambda^S}_{\iota_t}[f].
\end{equation}

Denote this optimal value by:

\begin{equation}
Q_t(x_t, u_t) = \mathcal{R}^{\lambda^S}_{\iota_t}[c_t + V_{t+1}] = -\frac{1}{\lambda^S} \log \mathbb{E}_{\iota_t}\left[ e^{-\lambda^S (c_t + V_{t+1})} \right].
\end{equation}

Moreover, the optimal transition kernel $\tau_t^*$ achieving this extremum is given by the exponential tilting formula from \eqref{eq:generalized-optimal-dist}:

\begin{equation}
\tau_t^* = \iota_t \frac{e^{-\lambda^S (c_t + V_{t+1})}}{\mathbb{E}_{\iota_t}[e^{-\lambda^S (c_t + V_{t+1})}]} = \iota_t \frac{e^{-\lambda^S (c_t + V_{t+1})}}{e^{-\lambda^S Q_t}}
\end{equation}

\begin{table*}[ht!]
\centering
\footnotesize
\setlength{\tabcolsep}{2pt}
\setlength{\extrarowheight}{0pt}
\renewcommand{\arraystretch}{3.5}
\caption{Dynamic programming recursions for control formulations}
\label{tab:dp_recursions}
\begin{tabular}{|c|c|c|c|c|c|c|}
\hline
 & \textbf{SOC} & \textbf{SP-SOC} & \textbf{RSOC} & \textbf{SP-RSOC} & \textbf{DOC} & \textbf{SP-DOC} \\
\hline
$V_t$ & 
$\min\limits_{\pi_t} \mathbb{E}_{\pi_t}[Q_t]$ &
$\min\limits_{\pi_t} \mathbb{E}_{\pi_t}[Q_t]{\,+\,}\frac{1}{\lambda^P}\mathbb{D}^{\pi_t}_{\rho_t}$ &
$\min\limits_{\pi_t} \mathbb{E}_{\pi_t}[Q_t]$ &
$\min\limits_{\pi_t} \mathbb{E}_{\pi_t}[Q_t]{\,+\,}\frac{1}{\lambda^S}\mathbb{D}^{\pi_t}_{\rho_t}$ &
$\min\limits_{\pi_t} \mathbb{E}_{\pi_t}[Q_t]$ &
$\min\limits_{\pi_t} \mathbb{E}_{\pi_t}[Q_t]{\,+\,}\frac{1}{\lambda^P}\mathbb{D}^{\pi_t}_{\rho_t}$ \\
 & & $= \mathcal{R}^{\lambda^P}_{\rho_t}[Q_t]$ & & $= \mathcal{R}^{\lambda^S}_{\rho_t}[Q_t]$ & & $= \mathcal{R}^{\lambda^P}_{\rho_t}[Q_t]$ \\
\hline
$Q_t$ & 
$c_t{\,+\,}\mathbb{E}_{\iota_t}[V_{t+1}]$ &
$c_t{\,+\,}\mathbb{E}_{\iota_t}[V_{t+1}]$ &
$\operatorname*{opt}^{\lambda^S}\limits_{\tau_t} \mathbb{E}_{\tau_t}[c_t{\,+\,}V_{t+1}]{\,+\,}\frac{1}{\lambda^S}\mathbb{D}^{\tau_t}_{\iota_t}$ &
$\operatorname*{opt}^{\lambda^S}\limits_{\tau_t} \mathbb{E}_{\tau_t}[c_t{\,+\,}V_{t+1}]{\,+\,}\frac{1}{\lambda^S}\mathbb{D}^{\tau_t}_{\iota_t}$ &
$c_t{\,+\,}V_{t+1}$ &
$c_t{\,+\,}V_{t+1}$ \\
 & & & $= \mathcal{R}^{\lambda^S}_{\iota_t}[c_t{\,+\,}V_{t+1}]$ & $= \mathcal{R}^{\lambda^S}_{\iota_t}[c_t{\,+\,}V_{t+1}]$ & & \\
\hline
$\pi^*_t$ & 
$\arg\min\limits_{\pi_t} \mathbb{E}_{\pi_t}[Q_t]$ &
$\rho_t \dfrac{e^{-\lambda^P Q_t}}{e^{-\lambda^P V_t}}$ &
$\arg\min\limits_{\pi_t} \mathbb{E}_{\pi_t}[Q_t]$ &
$\rho_t \dfrac{e^{-\lambda^S Q_t}}{e^{-\lambda^S V_t}}$ &
$\arg\min\limits_{\pi_t} \mathbb{E}_{\pi_t}[Q_t]$ &
$\rho_t \dfrac{e^{-\lambda^P Q_t}}{e^{-\lambda^P V_t}}$ \\
\hline
$\tau^*_t$ & 
$:= \iota_t$ &
$:= \iota_t$ &
$\iota_t \dfrac{e^{-\lambda^S (c_t{\,+\,}V_{t+1})}}{e^{-\lambda^S Q_t}}$ &
$\iota_t \dfrac{e^{-\lambda^S (c_t{\,+\,}V_{t+1})}}{e^{-\lambda^S Q_t}}$ &
$\iota_t$ &
$\iota_t$ \\
\hline
\end{tabular}
\end{table*}

Now the expression for $V_t$ simplifies to:

\begin{equation}
V_t(x_t) = \min_{\pi_t} \left[ \mathbb{E}_{\pi_t}[Q_t] + \frac{1}{\lambda^P} \mathbb{D}^{\pi_t}_{\rho_t} \right].
\end{equation}

Applying the dual representation \eqref{eq:generalized-risk-operator}:

\begin{equation}
V_t(x_t) = \mathcal{R}^{\lambda^P}_{\rho_t}[Q_t] = -\frac{1}{\lambda^P} \log \mathbb{E}_{\rho_t}\left[ e^{-\lambda^P Q_t} \right].
\end{equation}

The optimal policy achieving this minimum is given by the exponential tilting formula:

\begin{equation}
\pi_t^* = \rho_t \frac{e^{-\lambda^P Q_t}}{\mathbb{E}_{\rho_t}[e^{-\lambda^P Q_t}]} = \rho_t \frac{e^{-\lambda^P Q_t}}{e^{-\lambda^P V_t}}.
\end{equation}
% \newpage
We therefore have:

\begin{subequations}\label{eq:proof_recursions}
\begin{align}
V_t &= \min_{\pi_t} \left[ \mathbb{E}_{\pi_t}[Q_t] + \frac{1}{\lambda^P} \mathbb{D}^{\pi_t}_{\rho_t} \right] = \mathcal{R}^{\lambda^P}_{\rho_t}[Q_t], \\
Q_t &= \operatorname*{opt}^{\lambda^S}_{\tau_t} \left\{ \mathbb{E}_{\tau_t}[c_t + V_{t+1}] + \frac{1}{\lambda^S} \mathbb{D}^{\tau_t}_{\iota_t} \right\} = \mathcal{R}^{\lambda^S}_{\iota_t}[c_t + V_{t+1}], \\
\pi_t^* &= \rho_t \frac{e^{-\lambda^P Q_t}}{e^{-\lambda^P V_t}}, \\
\tau_t^* &= \iota_t \frac{e^{-\lambda^S (c_t + V_{t+1})}}{e^{-\lambda^S Q_t}}.
\end{align}
\end{subequations}

These are exactly the recursions in \eqref{eq:central_recursion}.

We have therefore shown:
\begin{enumerate}
    \item A terminal condition: $V_T = c_T$.
    \item If $V_{t+1}$ satisfies the recursions \eqref{eq:central_recursion}, then $V_t$ also satisfies them.
    \item The recursions \eqref{eq:central_recursion} hold for $V_{T-1}$.
\end{enumerate}

Therefore, by backward induction, the optimal solution of the central KL-regularized optimal control problem \eqref{eq:central_problem_full} is given by the recursions in Theorem \ref{thm:central_kl_problem}.
% \vspace{-15pt}

\section{Proof of Theorem \label{proof:compsitionality}}

The proof proceeds by induction. The base case $t=T$ holds by definition. Assume the decomposition holds at time $t+1$, so $z_{t+1} = \sum_{n=1}^N z_{t+1}^{(n)}$. Then:
\begin{align*}
z_t &= \mathbb{E}_{\rho_t}\!\left[r_t\,\mathbb{E}_{\iota_t}\!\left[z_{t+1}\right]\right] \\
&= \mathbb{E}_{\rho_t}\!\left[r_t\,\mathbb{E}_{\iota_t}\!\left[\sum_{n=1}^N z_{t+1}^{(n)}\right]\right] \\
&= \sum_{n=1}^N \mathbb{E}_{\rho_t}\!\left[r_t\,\mathbb{E}_{\iota_t}\!\left[z_{t+1}^{(n)}\right]\right] \\
&= \sum_{n=1}^N z_t^{(n)}.
\end{align*}

For the policy decomposition we get:
\begin{align*}
\pi_t^* &= \rho_t\,\frac{r_t\,\mathbb{E}_{\iota_t}[z_{t+1}]}{z_t} \\
&= \rho_t\,\frac{r_t\,\mathbb{E}_{\iota_t}\!\left[\sum_{n=1}^N z_{t+1}^{(n)}\right]}{\sum_{m=1}^N z_t^{(m)}} \\
&= \sum_{n=1}^N \frac{z_t^{(n)}}{\sum_{m=1}^N z_t^{(m)}} \cdot \rho_t\,\frac{r_t\,\mathbb{E}_{\iota_t}[z_{t+1}^{(n)}]}{z_t^{(n)}} \\
&= \sum_{n=1}^N w_t^{(n)} \pi_t^{(n)},
\end{align*}
where $w_t^{(n)} = \frac{z_t^{(n)}}{z_t}$ and $\sum_{n=1}^N w_t^{(n)} = 1$ since $z_t = \sum_{n=1}^N z_t^{(n)}$.

\end{document}